# Damage evaluation in a concrete gravity dam using Smoothed particle hydrodynamics

Tapan Jana[1], Amit Shaw[1], and L.S. Ramachandra[1]

[1] Department of Civil Engineering, Indian Institute of Technology Kharagpur, India

**Abstract.** An unexpected failure of a concrete gravity dam may cause unimaginable human suffering and massive economic losses. An earthquake is the main factor contributing to the concrete gravity dam's failure. In recent years, there has been a rise in efforts globally to make dams safe under dynamic loading. Numerical modeling of dams under earthquake loading yields substantial insights into dams' fracture and damage progression. In the present work, a particle-based computational framework is developed to investigate the failure of the Koyna dam, a concrete gravity dam in India exposed to dynamic loading. The dam-foundation system is considered here. The numerically obtained crack results in the concrete dam are compared with the available experimental results. The findings are consistent with one another.

## 1    Introduction

A concrete gravity dam plays a pivotal role in meeting the diverse demands of agriculture and industry. The structural integrity of such dams can be compromised by a multitude of factors, including seismic activity, thermal stress, differential foundation settlement, pre-existing structural flaws, and other variables. Among these factors, seismic-induced ground acceleration emerges as the primary catalyst for dam failure. In light of this, there has been a notable upswing in endeavors aimed at enhancing dam safety in recent years. The construction of dams in regions prone to seismic events necessitates meticulous planning and thorough preparatory measures. The evaluation of the dynamic response of the dam has been explored by many researchers employing the finite element method (FEM) ([1], [2], [3]) and extended finite element technique (XFEM) ([4], [5], [6]). Additionally, a prototype of a non-overflow monolith was experimentally investigated by Mridha et al. [7]. However, accurately determining free surfaces, deformable boundaries, and moving interfaces within the Eulerian framework presents a formidable challenge. In recent years, smoothed particle hydrodynamics (SPH) has gained much attention as a prominent numerical technique for evaluating concrete fractures due to its inherent property to take care of large deformation and material segregation. Das et al. [8] investigated the response of the Koyna dam to the seismic loading of the 1967 Koyna earthquake. The authors applied the continuum damage model proposed by Grady-Kipp. The authors mainly concentrated on the parametric study, i.e., the effect of excitation frequency, excitation amplitude, and the dissipation energy of the dam. The present work aims to test how effectively an SPH-based framework could identify dynamic load-induced crack patterns in a concrete dam.

## 2    SPH governing equations

Among the mesh-free approaches, SPH has found significant use in various domains. It was first developed by Lucy [9] and Gingold [10] to simulate astrodynamical problems. After that, it has been used in different problems, such as impact mechanics ([11], [12]), fracture mechanics ([13], [14], [15]), etc. The conservation equations for mass, momentum, and energy, based on the Lagrangian description, are represented in Einstein's indicial notations as,

$$\frac{d\rho}{dt} = -\rho \frac{\partial v^\beta}{\partial x^\beta} \tag{1}$$

$$\frac{dv^\alpha}{dt} = \frac{1}{\rho} \frac{\partial \sigma^{\alpha\beta}}{\partial x^\beta} \tag{2}$$

$$\frac{de}{dt} = \frac{\sigma^{\alpha\beta}}{\rho} \frac{\partial v^\beta}{\partial x^\beta} \tag{3}$$

where $\rho$ represents the density, $t$ represents the time, $x^\beta$ signifies the $\beta$-th element of the position vector, $v^\alpha$ represents the $\alpha$-th element of the velocity vector, $\sigma^{\alpha\beta}$ signifies the $\alpha, \beta$-th components of the Cauchy stress tensor, and $e$ stands for specific energy. In the context of SPH, the entire computational domain is represented by a collection of particles accomplished by applying any particle-based discretization technique. Let us say $N$ particles populate the discretized domain. The conservation equations in their discrete form can be expressed as

$$\frac{d\rho_i}{dt} = \sum_{j \in N} m_j (v_i^\beta - v_j^\beta) W_{ij,\beta} \tag{4}$$

$$\frac{dv_i^\alpha}{dt} = \sum_{j \in N} m_j \left( \frac{\sigma_i^{\alpha\beta}}{\rho_i^2} + \frac{\sigma_j^{\alpha\beta}}{\rho_j^2} - \pi_{ij} \delta^{\alpha\beta} \right) W_{ij,\beta} \tag{5}$$

$$\frac{de_i}{dt} = -\frac{1}{2} \sum_{j \in N} m_j (v_i^\beta - v_j^\beta) \left( \frac{\sigma_i^{\alpha\beta}}{\rho_i^2} + \frac{\sigma_j^{\alpha\beta}}{\rho_j^2} - \pi_{ij} \delta^{\alpha\beta} \right) W_{ij,\beta} \tag{6}$$

where $W_{ij,\beta} = \partial W(x_i - x_j, h)/\partial x_i^\beta$ is the gradient of the kernel function with smoothing length $h$. Here, a bell-shaped cubic B-spline function is used as the smoothing or kernel function. In the above equations, $\pi_{ij}$ is the artificial viscosity [19].

## 3   Fracture model and Pseudo-spring algorithm

In order to accurately capture the progressive damage within a material, it becomes imperative to integrate an appropriate fracture model into the computational framework. In this study, we adopt a fracture model based on the cohesive zone model, as devised by Kurumatani et al. [16], to simulate the fracture behavior of concrete. So, to incorporate the effects of material damage, one can write Hooke's law as $\sigma = (1 - D(\epsilon))E\epsilon$. Here, $D$ represents the damage variable, which provides a rough approximation of the damage sustained by a particular material. The damage variable ($D$) may be expressed as

$$D(\epsilon) = 1 - \frac{\epsilon_0}{\epsilon}\exp(-\frac{E\epsilon_0 h}{G_f}(\epsilon - \epsilon_0)) \tag{7}$$

where $\epsilon \geq 0$ is the maximum strain the material will undergo before failure occurs. In the context of damage evaluation, one now needs to capture the geometrical discontinuity in the material that may arise due to arbitrarily evolving and interacting cracks and moves through the material. In this study, the pseudo-spring technique ([13], [15], [17]) is utilized in conjunction with SPH to track damage in the material. Particles in pseudo-spring SPH are assumed to be connected with virtual springs to quantify the interaction between particles, as shown in Fig. 1. The extent of contact between any two particles is determined by the amount of damage accumulated in the connecting spring of those particles. When deterioration occurs, the interaction between particles weakens. When the attached spring completely fails, the interaction becomes zero.

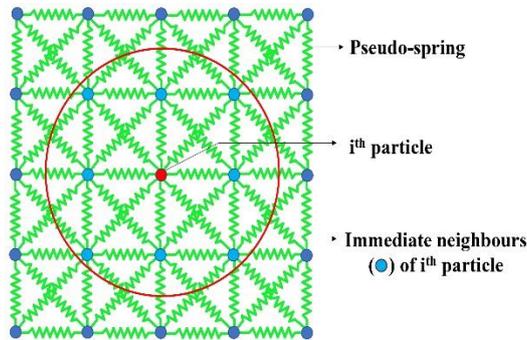

**Fig. 1.** Pseudo-spring connection

## 4   Dam configuration

The Koyna dam in Maharashtra, India, damaged by an earthquake in 1967, serves as the case study for our investigation. The dam is 103 m tall, 70 m wide at its base, and 14.8 m wide at its peak. A vertical cross-section of the Koyna dam-foundation system, as described in [1], is shown in Fig. 2(a). Fig. 2(b) shows the horizontal component of the Koyna earthquake acceleration data, which is normal to the dam's longitudinal axis, and Fig. 2(c) shows the vertical component. The coordinate value from these acceleration plots multiplied by gravitational acceleration (*g*) equals the applied acceleration value at the base of the dam for the damage estimation.

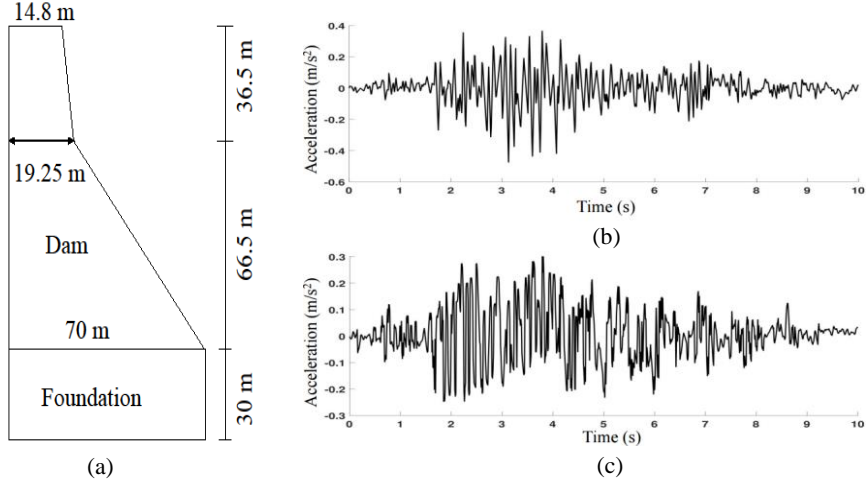

**Fig. 2.** (a) Koyna dam-foundation system (b) Koyna earthquake data (Horizontal) (c) Koyna earthquake data (Vertical)

## 5 Results for estimation of damage

The Koyna dam, in conjunction with its foundation, is considered in the present work. Specifically, our investigation involves modeling the dam's response under two distinct loading conditions: sinusoidal excitation and the seismic impact of the Koyna earthquake. Considering the particle spacing of 0.5 m, the total no of particles in the dam and foundation are 14521 and 8601, respectively. This precise particle arrangement facilitates a comprehensive examination of the structural behavior under the prescribed loading scenarios, enhancing the accuracy and reliability of our analysis. The smoothing length considered is 0.45m. Artificial viscosity coefficients ($\eta_1, \eta_2$) are 1 and 2. Rayleigh damping parameters ($\alpha, \beta$) are calculated as 1.616 and 0.0008, considering the first four fundamental frequencies from the work of [1]. The time step in compliance with the CFL criteria is $5 \times 10^{-6}$s. The total simulation time is 10s. The material properties and concrete fracture model parameters are considered from the work of [1] and [16], respectively. Young's modulus is 31.03 GPa, Poisson's ratio is 0.2, compressive strength is 31.9 MPa, tensile strength is 3.19 MPa, and Griffith fracture energy ($G_f$) is 100 N-m. The value of $\epsilon_0$ is 0.0001.

### 5.1 Damage under sinusoidal loading

A sinusoidal load of amplitude $0.1g$ and time period 0.3s is exerted on the dam's base in the horizontal direction. This load generates stress within the concrete material, leading to the initiation of damage in areas where stress levels exceed the concrete's tensile strength. Crack initiation initially takes place in the neck region due to pronounced stress concentration. Subsequently, crack propagation begins in the base region, albeit over a shorter distance. The crack follows the paths that are normal to the maximum principal stress direction and parallel to the direction of minimum principal stress. The contour plot of the dam's maximum and minimum principal stress at 10s is shown in Fig. 3. The simulated damage result is used for a comparative analysis, juxtaposing our simulated damage outcomes with the experimental findings of Niwa et al. [18]. This comparison is presented in Fig. 4, highlighting a remarkable alignment between our predicted fracture pattern and the experimental observations.

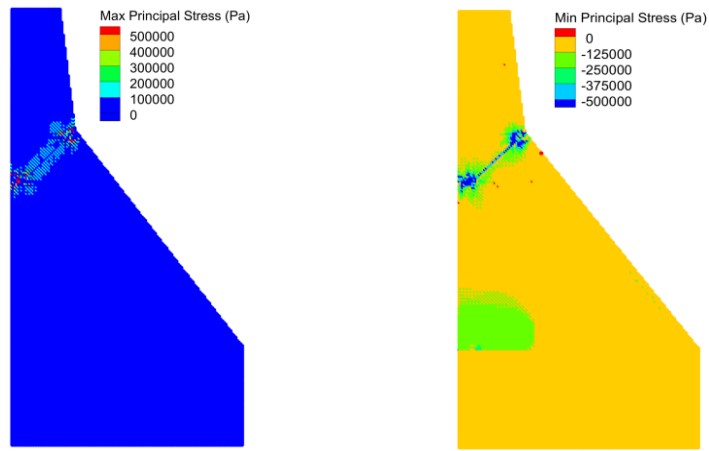

**Fig. 3.** Contour of maximum (left) and minimum (right) principal stress at 10 s

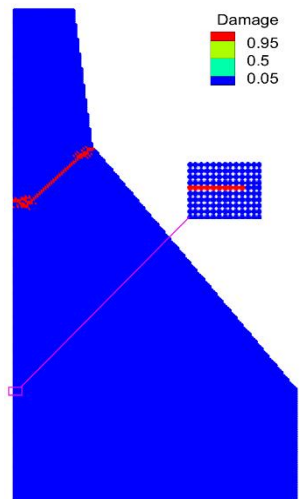 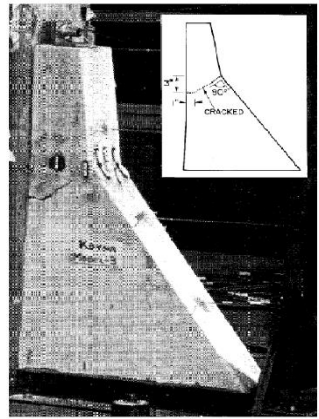

(a)  (b)

**Fig. 4.** Comparison of crack: (a) Simulated result (b) Experimental result from Niwa et al. [18]

### 5.2 Damage under earthquake loading

We now examine the performance of the Koyna dam when subjected to acceleration data recorded during the Koyna earthquake, as depicted in Fig. 2(a). The material characteristics and the numerical parameters are preserved from the last one. As anticipated, the neck is the first place for the crack initiation zone. Fig. 5 displays the contour plot depicting the maximum and minimum principal stress at the end of the simulation, i.e., 10s. Furthermore, Fig. 6 presents the crack profile at the conclusion of the simulation, i.e., 10s. The results are subsequently compared with the findings of Bhattacharjee and Leger [1] and Pekau et al. [3]. The observed crack profile at the neck and base matches well with the experiment findings. The appearance of a crack at the neck region aligns with the expectations due to pronounced stress concentration in that vicinity. Subsequently, a crack initiates at the base, following a straight horizontal trajectory for a distance of 10 m. Beyond this point, it takes on an inclined path, extending towards the downstream face at a height of 30m from the dam-foundation interface.

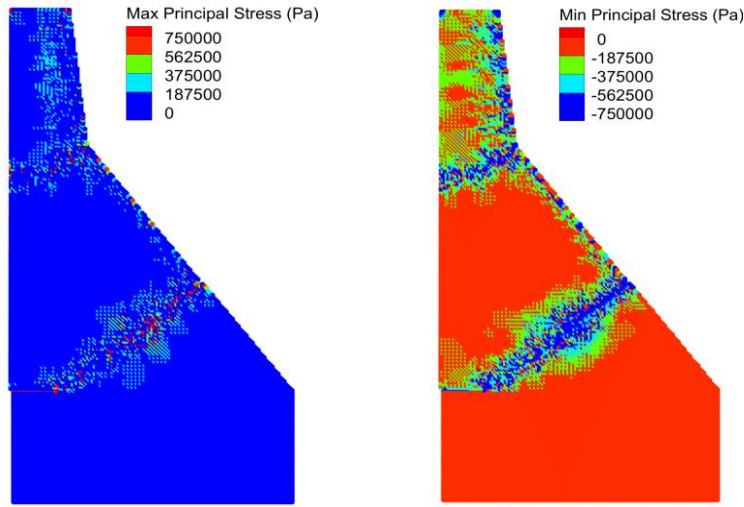

**Fig. 5.** Contour of maximum (left) and minimum (right) principal stress at 10 s

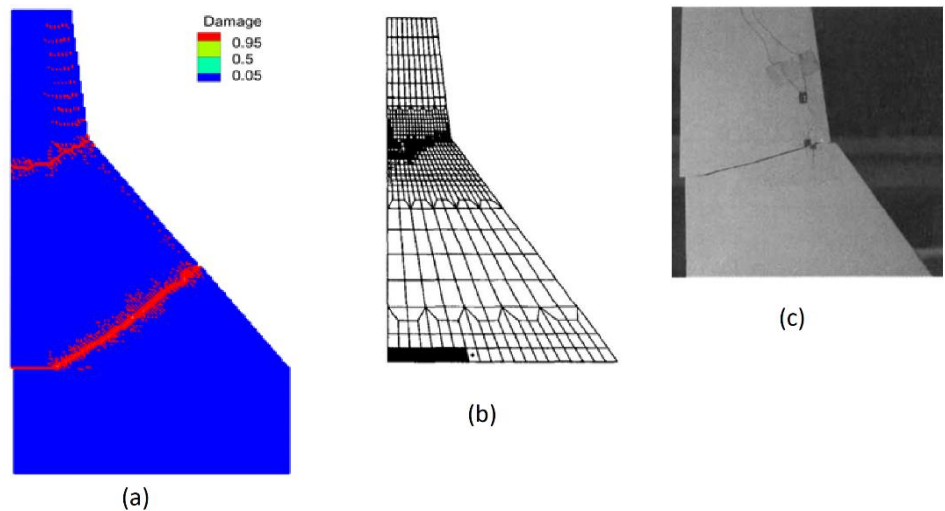

**Fig. 6.** Comparison of crack: (a) Simulated result (b) FEM result from Bhattacharjee and Leger [1] (c) Experimental result from Pekau et al. [3]

## 6   Conclusion

The present work investigates the evolution of cracks within a concrete gravity dam-foundation system, specifically focusing on the Koyna dam exposed to dynamic loadings. An SPH-based framework has been developed to simulate the Koyna dam under dynamic loadings. For the prediction of crack paths, a suitable concrete damage model has been incorporated. The simulated outcomes have been compared with available experimental results. The results demonstrate that the present developed framework can provide significant insights into evaluating the performance of a dam-foundation system exposed to any dynamic loads.